\documentstyle[12pt]{article}
\title{\large\bf CONTROLLING  STATIONARY FRONTS
IN TWO-DIMENSIONAL  REACTION-DIFFUSION SYSTEMS
}
\author{\normalsize \bf Yelena \,Smagina and Moshe Sheintuch  \\[2ex]
        \small {\it  Department of Chemical Engineering, Technion,} \\
        \small {\it Haifa, Israel,32000 }\\
        \small {\it email: cermssy@tx.technion.ac.il, Fax: 04-8230476}
       }
\date{}

\begin{document}
\maketitle

\pagestyle{empty}
\begin{abstract}
This paper considers new approach to control a stationary 
inhomogeneous  planar front solution of a nonlinear parabolic 
two-dimensional distributed (reaction-diffusion) 
system,  by using a gain point-sensor control with actuators that 
have the simplest possible  spatial dependence. The method is 
based on multivariable root-locus technique for the 
finite-dimensional approximation
of the original PDE model and  use the concepts of finite and
 infinite zeros of linear multidimensional
 system. 

\end{abstract}

                        \section{Introduction }

Nonlinear parabolic partial differential equations (PDEs), which typically describe reaction-diffusion 
systems, may admit spatially-dependent solutions like stationary fronts as well as spatiotemporal patterns.
 The latter can often be described as composed of slow-moving fronts, separated by domains of moderate changes. 
 This article is part of a research aimed to develop control theory for reaction-diffusion and
 reaction-convection-diffusion systems for which a certain patterned state is advantageous. 
Propagating fronts and patterned states may emerge in several technologies including 
catalytic reactors \cite{24},
 distillation processes \cite{10}, flame propagation and crystal growth \cite{1} (see also \cite{23} 
 for references)
 as well 
as in physiological systems like the heart \cite{18}.  Our interest lies in catalytic reactors, in
which stationary
 or moving fronts and spatiotemporal patterns have been observed and simulated in various systems like 
flow through a catalyst \cite{15}, fixed-bed reactors \cite{3,25}, reactors with flow reversal \cite{12} 
and loop 
reactors \cite{29}. The instabilities emerge due to the thermal effects in exothermic reactions, due 
to self-inhibition by a reactant and due to slow reversible modifications of the surface. 
The construction of a controller that enables to stabilize some inhomogeneous solutions in 
one-dimensional (1-D) 
reaction-diffusion and reaction-convection-diffusion systems is currently a subject of intensive 
investigation
 \cite{4,5,6,16,17,21,22}.  Yet, most catalytic reactors, as well as physiological systems like the heart, 
 exhibit a behaviour that can be described by two- or even three-dimensional models.
 
In the present work we are interested in stabilizing stationary (planar) fronts in 
a rectangular two-dimensional (2-D) domain in which a chemical reaction-diffusion process occurs. 
Previous studies of 1-D systems demonstrated that the simplest approach is by applying the point-sensor
control, in which a single space-independent actuator responds to a sensor that is located at the front
position \cite{21}. In the 2-D problem the point-sensor control with a single space-independent actuator 
cannot stabilize the front in a wide system. 
In that case it is necessary  to use control with several point sensors and actuators (space-independent and 
space-dependent). 
Here we apply the extension of the point-sensor control to multivariable case 
 by explore the  root-locus technique for multivariable system \cite{11}. The main advantage of this  control
is its insensitivity  to parameter uncertainties (robustness)  and  also small number of evaluated 
gain coefficients.  
 We apply the root-locus method to determine the minimal number of actuators ,
 their  spatial form and the location of the 
corresponding sensors that will assure the linear stability of the planar  front. The method uses 
the concepts of finite \cite{20,27} and infinite zeros \cite{11} of a linear multidimensional system; 
these notions are briefly explained below. The proposed approach is best suited for a systematic 
computer-aided search of the regulator form.

\section{Statement}

Consider the reaction-diffusion problem, in the $(z,r)$ rectangular domain  of length $L$ and width $R$ ,
 which is described by  a pair of coupled nonlinear parabolic PDEs 
\begin{equation}
\label{1a}
 y_t -y_{zz}-y_{rr} = P(y,\theta)+\lambda,\;\;\; \theta_t = \epsilon Q(y,\theta)
\end{equation}
 
subject to  no-flux boundary conditions :

\begin{equation}
\label{1b}
y_z(0,r)=0,\; y_z(L,r)=0,\; y_r(z,0)=0,\; y_r(z,R)=0
 \end{equation}
where the variable $y = y(z,r,t)$ typically represents the activator, $\theta = \theta(z,r,t)$ is the 
slow variable (localized inhibitor), $\lambda$ is a control variable that is introduced additively; 
$\epsilon < 1$ is the ratio of time scales.
 
We use the polynomial source functions
\begin{equation}
\label{2}
P(y,\theta)  = - y^3 +y+ \theta,\;\;\; Q(y,\theta) = -\gamma y  - \theta
\end{equation}
since they  adequately simulate the phenomenon of
 multiple steady states (bistable kinetics)
 due to thermal and autocatalytic effects which, in many chemical systems, induce the instabilities 
\cite{14} and since several analytical results are available \cite{13}. The instability 
stems  from the anticlinal arrangement of  $y$ and $\theta$ \cite{23}. 
In Eqn.(\ref{2}) $\gamma > 0 $ is a constant parameter.  
 System (\ref{1a}-\ref{2}) has often 
been used for studying pattern-formation in chemically reacting systems \cite{24}; similar systems were
 employed to describe pattern formation and control in cardiac systems.
 
	In this work we address the problem of stabilizing the planar front, in a 2-D reaction-diffusion system 
(\ref{1a}-\ref{2}),  at the position which in an open-loop system is stationary but unstable (i.e., 
at the middle, $z = L/2$ of the domain). Obviously, narrow systems behave like a one-dimensional system 
(Eqn. (\ref{1a}-\ref{2}), $y_{rr} = 0$) and admit a stationary front solution that is unstable and 
 typically oscillates or travels out of the system (see details in \cite{22}). In a sufficiently wide 2-D system,
 the front-line may
 undergo symmetry breaking so that in part of it the upper state expands while in other parts the 
lower state propagates. In that case we suggest to use control with several actuators located along 
the front at the points $(z,r)=(L/2,r_d),d=1,\ldots,\eta $   and apply a general feedback control law of the form

\begin{equation}
\label{3}
 \lambda =  k\sum_{d=1}^{\nu} [y(L/2,r_d,t)-y^*_d]\psi_d (z,r)
\end{equation}
where $k<0$ is a scalar gain coefficient,  $ y(L/2,r_d,t)-y^*_d$  are deviations of 
the sensors from the set points; $ y^*_d =y^*(L/2,r_d)$; $ \psi_d (z,r)$ 
 are  some space-dependent  functions  that may imitate  the eigenfunctions. 
We will seek control (\ref{3}) with the simplest space-independent or space-dependent actuator 
functions $ \psi_d (z,r)$   and  minimal number  of sensors.

\section{Linear analysis}

For design of control (\ref{3}) we use a linearized truncated version of  PDEs (\ref{1a}-\ref{2}).  Linearization of 
(\ref{1a}-\ref{2}) for  $\lambda = 0$ around the steady state solution  $y_s =  y_s (z,r)$, $\theta _s =\theta_s(z,r)$
yields 
\begin{equation}
\label{4}
\bar{y}_t - \bar{y}_{zz} -\bar{y}_{rr}  = (-3\bar{y}^2_s +1)\bar{y} + \bar{\theta} + \lambda, \qquad
 \bar{\theta}_t = -\epsilon \gamma \bar{y} -\epsilon \bar{\theta}
\end{equation}
which is lumped by the Galerkin method:  We expand the deviations $\bar{y}= \bar{y}(z,r,t)$,
$\bar{\theta}= \bar{\theta}(z,r,t)$ and $\bar{y}(z^*,r_d,t)= y(z^*,r_d,t)-y_s(z^*,r_d)$ as
\begin{equation}
\label{5}
 \bar{y} (z,t) =  \sum_{e} a_e (t)\phi_e (z,r), \;
 \bar{\theta} (z,t) =  \sum_{e} b_e (t)\phi_e (z,r), \;
 \bar{y}(z^*,r_d,t)= \sum_{e} a_e (t)\phi_e (z^*,r_d)
\end{equation}
where $z^*=L/2$ and the functions $ \phi_e (z,r)$ are the eigenfunctions of the problem
\begin{equation}
\label{6}
\phi_{zz} (z,r)+ \phi_{rr} (z,r) = -\lambda \phi (z,r),\;\;\;
\phi_z (0,r)=\phi_z(L,r)= \phi_r (z,0)=\phi_r (z,R)=0
\end{equation}
with the corresponding eigenvalues
\begin{equation}
\label{7a}
\lambda_{e(ij)}=\pi ^2[(i-1)^2/L^2+(j-1)^2/R^2]
 \end{equation}
and  the eigenfunctions
\begin{equation}
\label{7b}
\phi_{e(ij)}=\frac{\rho}{\sqrt{LR}}cos\frac{(i-1)\pi z}{L}cos\frac{(j-1)\pi z}{R},\;\;\; e=e(ij)=1,2,\ldots
\end{equation}
($\rho=1$ when $i=j=1$, $\rho=2$ when $i,j>1$ and $\rho=\sqrt{2}$ when $i=1$, $j>1$  or $j=1$,$i>1$,
 see \cite{23}  for derivation ). Substituting (\ref{5}) into (\ref{4}),(\ref{3}) with  set points
 $y_d^* =y_s(L/2,r_d)$
 and integrating with a weight eigenfunctions $\phi_e(z,r)$ results in the spectral representation of
 closed-loop linearized system (\ref{4}),(\ref{3}) 
\begin{equation}
\label{8}
\dot{a}_e = -\lambda _ea_e + \sum_{m}J_{em}a_m + b_e + 
k\sum_{d=1}^{\eta}\int_{0}^{L}\int_{0}^{R}\psi_d (z,r)\sum_{f}a_f\phi_{f(kl)}(L/2,r_d)\phi_{e(kl)}(z,r)dzdr
\end{equation}
\begin{equation}
\label{9}
\dot{b}_e = \epsilon(\gamma a_e + b_e)
\end{equation}
where
\begin{equation}
\label{10}
J_{em} = \int_{0}^{L}\int_{0}^{R}(-3y_s^2 + 1)\phi_{m(ij)}\phi_{e(kl)}dzdr\;\;\; e,m=1,\ldots
\end{equation}
Denoting
\begin{equation}
\label{11}
h_{df} = \phi_f (L/2,r_d)
\end{equation}
\begin{equation}
\label{12}
\beta_{ed} = \int_{0}^{L}\int_{0}^{R}\psi_d (z,r)\phi_{e(kl)}(z,r)dzdr
\end{equation}
we can  rewrite (\ref{8})   as follows
\begin{equation}
\label{13}
\dot{a}_e = -\lambda _ea_e + \sum_{m}J_{em}a_m + b_e + 
k\sum_{d=1}^{\eta}\beta_{ed}\sum_{f}a_f h_{df},\;\;\; e,f=1,2,\ldots
\end{equation}
This system  (Eqn. \ref{13}, \ref{9})  may be presented  in the usual vector-matrix form  as the linear
 infinite-dimensional dynamical  system  with $\eta$ - dimensional  input  $v$ and  output $w$ vectors

\begin{equation}
\label{14} 
\left[ \begin{array}{c} \dot{a}\\ \dot{b} \end{array} \right]   \;=\;
\left[ \begin{array}{cr} -\Lambda +J &  I \\ -\epsilon \gamma I & -\epsilon I
  \end{array} \right] \left[ \begin{array}{c} a\\ b \end{array} \right] +
       \left[ \begin{array}{c}
\beta \\ O \end{array}  \right] v
\end{equation}
\begin{equation}
\label{15} 
w=Ha
\end{equation}
closed by  the finite-dimensional output feedback 
\begin{equation}
\label{16} 
v=kI_{\eta}w
\end{equation}
where  $a(t) =[a_e], b(t) = [b_e]$  are the infinite dimensional vectors ($e=1,2,\ldots$);  $v$ and $w$
 are finite-dimensional $\eta$  vectors,  the matrix $\beta =[\beta_{ed}]$ has $\eta$ infinite-dimensional columns 
($e=1,2,\ldots; d=1,\ldots,\eta$)  and the matrix $H=[h_{df}]$  has $\eta$ infinite-dimensional rows 
 ($d=1,\ldots,\eta; f=1,2,\ldots$); $I_{\eta}$ is unity $\eta \times \eta $  matrix; 
$\Lambda =diag(\lambda_1,\lambda_2 , \ldots )$, $ I = diag(1, 1, \ldots )$ , $J = [J_{em} ],  e,m=1,2,\ldots $ 
 are infinite-dimension matrices and  $k$  is scalar gain coefficient. 

Let us study the correlation between  the form of original control (\ref{3}) and the input and output  structures
 of  Eqns.(\ref{14}-\ref{15}): If all $r_d,d=1,2,\ldots,\eta$  in (\ref{3}) have different values, 
 then the parameter $\eta$  assigns the dimension of input and output vectors of  Eqns. (\ref{14}), (\ref{15}) 
and the sensor positions stipulate the structure of the  output  matrix  $H$ (see Eqn.\ref{11}). Besides, 
the form of the actuator functions $\psi_d(z,r)$ influences the matrix $\beta$ (see Eqn \ref{12}).
 Hence, in view of the original statement the problem can be stated as follows:
 
{\bf{Problem 1}}. For the linearized infinite-dimensional ODEs system (\ref{14}-\ref{15}) it is necessary to 
find the matrices $\beta$ and $H$  with the minimal number of columns and rows, respectively, such that 
 the finite-dimensional output feedback control (\ref{16}) stabilizes the closed-loop system.
 
The obvious way of designing of a finite-dimensional control (\ref{16}) for an infinite-dimensional system 
(\ref{14}-\ref{15}) is to use a truncated (finite-dimensional) approximation of the PDEs. To truncate the system we
 capitalize on  the dissipative nature of the parabolic PDEs: the truncation order $N$ is estimated by calculating 
the leading eigenvalues of the dynamics matrix  (\ref{14}) (see \cite{23} for details). As follows from \cite{2}, 
for a sufficiently large truncated order $N$,  with actuator functions $\psi_d(z,r)$  that coincide with 
 eigenfunctions $\phi_e (z,r)$, the finite-dimensional control guarantees the stability of infinite-dimensional
 system.  Below we imply that system (\ref{14}-\ref{15}) is  the  finite-dimensional analog of the original
 (\ref{14}-\ref{15}) with truncated order $N$.

\section {Root-locus control design}

Design of control (\ref{3}) implies the determination of the minimal number of actuators to be employed, 
their spatial form and the location of the corresponding sensors that will assure the linear stability 
of the planar front.  We can  manipulate the number and r-coordinate positions ($r_1, r_2,
\ldots,r_{\eta}$  ) of the sensors and
 the form of the  actuator functions $\psi_d (z,r)$. The former affects  the matrix $H$ (Eqn. \ref{11})
 and the latter  influences  the matrix  $\beta$ (Eqn. \ref{12}).
  
For simplification of the search for a  matrix $\beta$  structure we assign the eigenfunctions (\ref{7b})
 as actuator functions, i.e. $\psi_d(z,r) \sim \phi_e (z,r), e=1,2,\ldots $. Thus from the relation  between 
the form of the actuator functions $\psi_d(z,r)$ 
and the structure of the  matrix $\beta$(see Eqn.\ref{12}) it follows that every  $d$-th  column of the matrix
 $\beta$  will contain only a single non-zero element $\beta_{ed}$ . Moreover, such a structure of $\beta$ 
 satisfies  the 
above-mentioned  Balas's restrictions \cite{2} that ensure the validity of control of an infinite-dimensional
 system by a finite- dimensional controller. 

Therefore, we  propose the following steps for design of control (\ref{3}): (i)  Assign sensor  positions along
 the front $(z,r)=(L/2,r_d), \; d=1,\ldots \eta \;\;$  and calculate the matrix $H$.  (ii) Seek the matrix $\beta$
  and gain $k$  that assures linear stability  
of the closed-loop system (Eqns.\ref{14}-\ref{16}). (iii) Finally, find the eigenfunctions $\phi_e(z,r)$
 that corresponds to the    matrix  $\beta$  obtained. 
  
For fulfilling step (ii)  we use a root-locus technique \cite{11}  which is based  on an analysis  of the 
finite zeros  and infinite zeros of the open-loop system (\ref{14}-\ref{15}). This method  uses the following
 known property of closed-loop linear system with feedback (\ref{16}) (see \cite{9}): as the feedback 
gain increases towards infinity a part of the closed-loop eigenvalues remain finite and approaches 
the values which are referred to as  finite system zeros \cite{27},\cite{20} 
(see also Appendix for definitions)
 while the remainder are located at the points at infinity and  are known as infinite zeros \cite{11}. 
 Therefore, we propose to seek a suitable matrix $\beta$  by the repeatedly calculating the finite zeros 
of open-loop ODEs (\ref{14}-\ref{15}), with different input matrices, and finding the one that ensures that  
the leading finite system zeros are negative. Then we assign a sufficienly large 
negative gain coefficient  $k$.\footnote{Here
 we must note that large perturbations of the front will cause the control variable  
 to exceed the bistability domain of $P = 0$ . So the gain value coefficient  $k$ must be restricted below 
the value ($k_o$)  that ensures stability of the closed-loop system.}

Finally,  Problem 1  may be reformulated as follows: 

{\bf{Problem 2.}} For the linearized truncated system (\ref{14}-\ref{15}) with assigned  output matrix $H$ 
it is necessary to find the matrix $\beta$  that has a single nonzero element in every column such 
that the above-mentioned system has leading finite system zeros  in the left -half of the complex plane
 (`negative` zeros)\footnote{Similarly, we denote  zeros in the right-half  of the complex plane  by 
 `positive` ones.}.

 Let us note that such input and output matrices  $H$ and $\beta$ result a minimum phase control system 
 because this system contains finite zeros in the left-half of the complex plane.
  
The following assertion is needed to avoid cases when Problem 2 has no solution .

{\bf{Assertion 1.}} For assiged matrix $H$ the problem has no solution for any matrix $\beta$  if and only if `positive`
  finite system zeros of  (\ref{14}), (\ref{15})  are  output-decoupling zeros \cite{19}.
 
The proof follows from definition of  decoupling zeros (see Appendix).

{\bf{Remark 1.}}  If the shape of actuator distribution functions $\psi_d (z,r), d=1,2,\ldots $
  are preassigned  from technical constraints (i.e. the matrix $\beta$ is given) then  we need  to find $\eta$
 sensor locations (i.e. the matrix $H$ ) which provide the 'negative' leading finite system zeros of  system 
(\ref{14}-\ref{15}). The solvability of this problem is formulated as follows.

{\bf{Assertion 2.}} For assiged matrix $\beta$ the problem has no solution for any matrix $H$  if and only if  `positive`  finite system zeros 
of  (\ref{14}), (\ref{15})  are  input-decoupling zeros \cite{19} .
 
{\bf{Remark 2.}} It is necessary to choose the sensor locations in $r$-direction so that the  infinitely increasing
 eigenvalues of  the closed-loop high gain system (\ref{14}-\ref{16}) (infinite zeros) tend  to infinity along 
asymptotics with a negative real angle. This condition is guaranteed if $det(H\beta \neq) 0$  and  all eigenvalues
 of the $\eta \times \eta $ matrix $H\beta$ are positive \cite{11}. If  $H\beta$  has several negative eigenvalues 
then we need to introduce a  nonsingular precompensator $M$  to  (\ref{16}) such that  the new control
\begin{equation}
\label{17} 
v=kI_{\eta}Mw
\end{equation}
ensures above  property for $MH\beta$ . Such operation does not changes the finite system zeros of 
(\ref{14}), (\ref{15})  because they are invariant to any nonsingular transformation of output \cite{27},\cite{20}.

Therefore, the general strategy of the method is as follows: At first we need to check that assigned sensor
 positions (or matrix $H$) ensure  that output-decoupling zeros of a pair $( A, H )$ are 'negative'  (see Assertion 1).
 Then we seek a suitable matrix $\beta$ by repeating calculations of the finite system zeros of  open-loop ODEs 
(\ref{14}), (\ref{15})  with different matrices $\beta$  and finding one that ensures that  the leading finite 
system zeros are negative and $det(H\beta)\neq 0$. If necessary, we find the precompensator  $M$  which rearranges 
the sensor position in $r$-direction  so that the infinite zeros tend to infinity along asymptotes with a negative 
real axis angle.

 To demonstrate this procedure we apply it for design of control (\ref{3})  for  PDEs (\ref{1a}-\ref{2}) of length $L$ 
 and various widths.

 We start by analyzing the effectiveness of the simplest control law, a single space-independent actuator 
 (Eqn. \ref{3}, $\eta=1$, $\psi_1(z,r)=1$)
\begin{equation}
\label{18} 
 \lambda(t)=ky((L/2,r_1,t)-y_s^*)
\end{equation}
where $y_s^*=y_s(L/2,r_1)$ coincides with the steady-state value of the problem .
 The spectral representation of  the closed-loop system (\ref{1a}-\ref{2}),(\ref{18}) 
 is a single-input, single-output  Eqns.(\ref{14}-\ref{15})  with a column  vector $\beta =[\beta_{11},0,\ldots]$ 
with $\beta_{11}=1/\sqrt{LR}$  and  
row  output vector $H$ = $h = [h_1,h_2,\ldots]$.
Thus here the shape of the actuator  is assigned. Then if  input-decoupling zeros of a pair $( A, \beta )$ are 'negative'  
 we need to find  a sensor position. At first consider a  sensor situated at the domain center 
($z^*=L/2, r_1=R/2$). The analysis of leading system zeros of the related  linearized truncated system obtained 
shows that they possess negative real parts for $0<R<R_{cr}$ for some value $R_{cr}$ and positive real parts 
 for $R>R_{cr}$. So, control with one space-independent actuator 
is effective only  for systems of width $R \leq R_{cr}$ (narrow systems) that have 'negative' zeros. 
Changing the sensor position in the $r$ -direction   alters 
the vector $h$  (Eq.\ref{11} ) and as a consequence may move the finite system zeros (for details see \cite{28})
and   $R_{cr}$. However, the position at the domain center assures  the maximal value of
$R_{cr}$.  

For wider systems ($R>R_{cr}$) it is necessary to use space-dependent actuators. We try  to apply  control in the form
\begin{equation}
\label{19}
 \lambda =k\sum_{d=1}^{2} [y(L/2,r_d,t)-y^*_s]\psi_d (z,r)
\end{equation}
with one space-independent actuator  ($\psi_1(z,r)=1$) and another space-dependent one, 
$\psi_2(z,r) \sim \phi_e(z,r)$  where  $\phi_e(z,r)$  is eigenfunction (\ref{7b}) 
chosen from the series eigenfunctions  ordered in an increasing
 order of the appropriate eigenvalues. Introducing two sensors at positions $(L/2,r_1)$ and $(L/2,r_2)$ 
 we  calculate the $2\times N$  matrix  $H$  by Eqn. (\ref{11}). 
If  output-decoupling zeros of a pair $( A, H )$ are 'negative'  then we need to evaluate the finite system zeros of 
systems (\ref{14}-\ref{15}) with above output matrix $H$ and different $N\times 2$ matrices $\beta =[\beta_{ed}]$,
$e=1,\ldots,N ,d=1,2$ with assigned elements of the first column
($\beta_{11}=1/\sqrt{LR},\;\; \beta_{e1}=0,\;e=2,\cdots,N$ )  and undetermined elements
 of the second column.

Thus it is necessary to find a single  nonzero element from the second column of $\beta$ which ensures that the leading
 finite zeros of system (\ref{14}-\ref{15}) are `negative` ones. If  such $\beta$ does not exist 
we need to change the matrix  $H$  by shifting  the positions of sensors in $r$-direction and begin  the search of
$\beta$  once more.  If an appropriate  matrix $H$ does not exist then control (\ref{19})  is not effective and 
it is necessary to increase the number of sensors (and actuators) to three and so on.
 
\section {Application}

 		Let us apply the above method for stabilization of planar front solution of system (\ref{1a}-\ref{2})
(with  $L=20$, $\gamma=0.45$,  $\epsilon=0.1$ ) of various widths. The leading eigenvalues of the truncated
 version (Eqn. \ref{14},  $ N = 23$) show two real unstable eigenvalues ( $0.35$ and $0$) for all $R$'s 
and two complex eigenvalues  with  real parts that becomes positive for $R \geq 5.6$( see Fig. 4a in \cite{23}).
 At first  we try to apply control (\ref{18})  with one space-independent actuator and 
the sensor at ($z^*=10$, $r_1=R/2$).
  The analysis of zeros of open-loop system (\ref{14}-\ref{15}) discovers 
 two leading (and complex) finite zeros with negative real parts for $0<R<5.6$ and positive real parts for
$R\geq5.6$ (see Fig.4b in \cite{23}). Hence, system  width $R_{cr}=5.5$  is a critical one for our  ability to do 
control with a space-independent actuator.

For wider systems ($R>5.5$) we apply the two actuator control  (\ref{19}). Assigning  two sensors in some 
positions ($z^*=10,r_1$),  ($z^*=10,r_2$) ($r_1 \neq r_2$) and  calculating the finite zeros  of 
 relevant two input, two-output systems (\ref{14}-\ref{15}) with different input  matrices   we find that when 
 $\beta_{24} \neq 0$ the leading finite zeros have negative real parts  ($-0.1024$). This  correspondents to 
the fourth ($e = 4$)  eigenfunction $\phi_4 \sim cos(\pi r/R)$  in series of the ordered eigenfunctions :
 The first six ordered eigenfunctions  are  $\phi_1 =1$, $\phi_2 \sim cos(\pi z/L)$, $\phi_3 \sim cos(2\pi z/L)$,
 $\phi_4 \sim cos(\pi r/R)$, $\phi_5 \sim cos(3\pi z/L)$, $\phi_6 \sim cos(2\pi z/L)cos(\pi r/R)$.
 Consequently control (\ref{19})  becomes
\begin{equation}
\label{21} 
 \lambda(t)=k\{[y(10,r_1,t)-y_s(10,r_1)] + [y(10,r_2,t)-y_s(10,r_2)]cos(\pi r/R \} 
\end{equation}
      
Then, it is necessary to choose $r_1$, $r_2$  in (\ref{21}) so that  Remark 2 is satisfied. For our case we 
need to use the sensor locations with $r_1 >r_2 $ which is equivalent to introducing the $2\times 2$   precompensator 
$ M = \left[ \begin{array}{cr} 0 & 1 \\ 1& 0 \end{array} \right]$. 

We verify the proposed methodology by simulating system (\ref{1a}-\ref{2}) using one space- independent actuator
 (\ref{18}) for narrow ($R<R_{cr}$) and two-actuator control  (\ref{21}) for wider system  (see \cite{23} for plots). 

\section { Conclusion remarks}

The stabilization of planar stationary fronts in a two-dimensional rectangular domain,
 in which a diffusion-reaction systems occurs, is studied using a two-variable PDEs model  
for which some analytical results are available. We consider the simplest control strategy based on sensors
 placed at the designed front line position and measure deviations from a local state, and actuators that 
are spatially-uniform or space dependent. We present a systematic control design that determines
 the number of required sensors and actuators, their position and their form.
 The control design is  corroborated by linear analysis of a lumped truncated model and concepts of finite 
and infinite zeros of linear multidimensional systems. The method is best suited for a systematic computer-aided
 search of the regulator form. 

 \begin{center} {\bf APPENDIX}\end{center}

 Consider a general linear multivariable finite-dimensional dynamic system described
 by the set of state-space equations
 $$ \dot{a}(t)=Aa(t)+Bv(t)\;\;\;\;\;\;    w(t)=Ha(t)$$
with  $a$ , $v$ and $w$ are the $n$ , $r$  and  $l$  dimensional state, input and output vectors  and 
 $A$, $B$ and $H$ are  constant matrices of appropriate dimensions. 
 
{\bf{Definition 1.}}  The  finite zeros (system zeros) of the  above system are determined  as  the set 
of complex $s_i$  for which the rank of the system matrix 
$$ P(s) = \left[ \begin{array}{cr} s_iI-A & -B \\ H & O \end{array} \right]$$ 
is  reduced.

{\bf{Definition 2.}} The input-decoupling zeros are defined  as the set of the complex variable $s_i$ 
at which the row rank of the matrix  $\left[ s_iI-A , -B  \right]$ is reduced.

{\bf{Definition 3.}} The output-decoupling zeros are defined  as the set of the complex variable $s_i$ 
at which the row rank of the matrix $\left[ s_iI-A^T , H^T   \right]$    is reduced.

{\bf{Evaluation of input-decoupling zeros.}} The set of input-decoupling zeros may be calculated from Definition 2.
 To avoid the operations with complex numbers we can apply an alternative method which uses the property:
 input-decoupling zeros coincide with the uncontrollable eigenvalues of the matrix $A$. Since the latter
 eigenvalues are invariant under a proportional state feedback: $v=Ka$ then they may be calculated as those 
eigenvalues of the closed-loop matrices $A+BK_j$   which are invariant with respect to any gain matrices $K_j$ 
with finite elements.  Let us analyze the leading eigenvalues of the $N$-truncated system (\ref{14}) acted 
by the state control $v=K_j a$ with a random $N$ row vector $K_j$. It is evident that leading eigenvalues, which are
 invariant with respect to this control, are input-decoupling zeros of (\ref{14}), (\ref{15}).

\begin{center}
{\bf ACKNOWLEDGMENT}
\end{center}

 	M.S. acknowledges the Minerva Center of  Nonlinear Dynamics for support.

\end{document}